\documentclass[11pt]{article}
\usepackage[T1]{fontenc}
\usepackage[utf8]{inputenc}
\usepackage{lmodern,amsmath,amssymb,amsthm,mathtools}
\usepackage[margin=1in]{geometry}
\usepackage{booktabs,microtype,enumitem,float}
\usepackage[colorlinks=true,linkcolor=blue!50!black,citecolor=blue!50!black,urlcolor=blue!50!black]{hyperref}
\usepackage{xcolor}
\hypersetup{pdftitle={Twelve normal forms for symmetric magic squares of order four},pdfauthor={Xinru Jiang and Yueming Zhong},
  pdfkeywords={symmetric magic squares, Hilbert basis, half-open cones, Ehrhart theory}}
\newtheorem{theorem}{Theorem}[section]
\newtheorem{lemma}[theorem]{Lemma}
\newtheorem{proposition}[theorem]{Proposition}
\newtheorem{corollary}[theorem]{Corollary}
\theoremstyle{remark}
\newtheorem{example}[theorem]{Example}
\newcommand{\N}{\mathbb Z_{\ge0}}
\newcommand{\Z}{\mathbb Z}
\newcommand{\R}{\mathbb R}
\newcommand{\cone}{\operatorname{cone}}
\newcommand{\conv}{\operatorname{conv}}
\newcommand{\Hilb}{\operatorname{Hilb}}
\newcommand{\intr}{\operatorname{int}}

\newcommand{\calS}{\mathcal S}
\newcommand{\calF}{\mathcal F}
\numberwithin{equation}{section}
\allowdisplaybreaks[2]
\title{Twelve normal forms for symmetric magic squares of order four}
\author{Xinru Jiang\thanks{Email: \texttt{chiang\_11110@163.com}.}\qquad
Yueming Zhong\thanks{Corresponding author. Email: \texttt{zhongyueming107@gmail.com}.}\\[4pt]
\small Key Laboratory of Engineering Modeling and Statistical Computation of Hainan Province\\
\small School of Mathematics and Statistics, Hainan University\\
\small Haikou 570228, P.R. China}
\date{}
\begin{document}
\maketitle
\begin{abstract}
We give twelve disjoint normal forms for symmetric $4\times4$ nonnegative
integer matrices whose rows, columns, main diagonal, and antidiagonal have a
common sum. Each class is a translate of a free submonoid on five generators.
A decision table and inverse formulas recover the class and its coefficients
from the entries. Every finite disjoint decomposition of the additive
monoid of these squares into translates of free submonoids requires at
least twelve classes, even if the ranks of the submonoids vary.
The construction gives a regular unimodular triangulation of the solution
cone and shows that its twelve primitive extreme-ray generators
form the Hilbert basis. It also yields a multivariate generating function
and a counting quasipolynomial of least period four in the magic sum.
\end{abstract}



\section{Introduction and main results}

Let $\calS$ be the additive monoid of symmetric matrices
$X=(x_{ij})\in\N^{4\times4}$ whose rows, columns, main diagonal, and
antidiagonal all sum to the same integer $N$. We call $N$ the
\emph{magic sum}. Entries may repeat, and the positions are labeled.

Beck, Cohen, Cuomo, and Gribelyuk studied this model using Ehrhart theory
\cite{BCCG03}. In Section~4 of their preprint, they report twelve vertices
and denominator four for the section of magic sum one
\cite{BCCGv2,BCCGv3}. Their displayed counting formula differs from the
one obtained here. We compare the formulas in
Section~\ref{subsec:comparison}, starting with a direct count at magic sum
two.

Stanley uses the term symmetric magic square for a nonnegative integer
matrix equal to its transpose with all row and column sums equal. For this
semimagic model, Carlitz gave the generating function in order four
\cite{Carlitz66,Stanley76}, and Ahmed gave the Hilbert basis of the solution
cone \cite[Section~4.4.1, Figure~4.9]{Ahmed04}. The two generating
functions are related by $F(t)=(1-t)^2G(t)$, where $G$ records the
common row sum in the semimagic model and $F$ records half the magic sum
in our model. A matching decomposition gives a direct explanation in
Section~\ref{subsec:countproof}. For the use of Hilbert
bases in related magic-square and magic-cube problems, see
\cite{Stanley73,ADH03}. Xin gave two normal forms for magic squares of
order three with distinct entries, up to rotations and reflections
\cite{Xin08}.

Stanley's theorem on homogeneous linear Diophantine systems already
gives a disjoint decomposition into translates of full-rank free
submonoids \cite[Theorem~5.2]{Stanley82}, following the polyhedral
construction in \cite[Section~1]{Stanley80}. For constructions based on
half-open simplicial cones and fundamental parallelepipeds, see
\cite[Section~5.3]{BIS16}. The same work gives a 48-class decomposition
for the rank-eight model of order four without transpose symmetry
\cite[Table~4]{BIS16}. Here four classes are selected by comparisons
among the entries and half the magic sum. In the remaining region, two
signs determine a change of coordinates, followed by one comparison in
a three-dimensional cone. These comparisons determine the normal form
and its coefficients. The poles of the counting series and
the degrees of the extreme-ray generators give the lower bound of
twelve classes.

We use five free entries, in the order $u=(b,c,d,f,g)$, and the linear map
\begin{equation}\label{eq:phi}
\Phi(u)=
\begin{pmatrix}
d+2f-b-c&b&c&d\\
b&2d+f-b-g&f&g\\
c&f&b+g-f&2d+2f-b-c-g\\
d&g&2d+2f-b-c-g&b+c-d
\end{pmatrix}.
\end{equation}
Its magic sum is $2(d+f)$. Let
\begin{equation}\label{eq:cone}
\begin{split}
C=\{(b,c,d,f,g)\in\R_{\ge0}^{5}:\;&b+c\ge d,\quad b+g\ge f,\\
&b+c\le d+2f,\quad b+g\le2d+f,\\
&b+c+g\le2d+2f\}.
\end{split}
\end{equation}
As shown in Section~\ref{subsec:coordinates}, $\Phi$ identifies $C\cap\Z^5$
with $\calS$. Set
\[
\deg u=d+f,\qquad \deg X:=\deg(\Phi^{-1}(X))\quad(X\in\calS).
\]
Thus the degree of a matrix is half its magic sum. Put $V_i=\Phi(v_i)$,
where the vectors $v_i$ are given in Table~\ref{tab:rays}; the corresponding
matrices appear in Appendix~\ref{app:matrices}. The section $P=C\cap\{\deg=1\}$ corresponds to magic sum two.
The map $\Phi$ sends $\tfrac12P$ to the section of magic sum one.

\begin{table}[H]
\centering
\caption{The twelve primitive generators.}\label{tab:rays}
\begin{tabular}{c@{\qquad}c@{\qquad}c|c@{\qquad}c@{\qquad}c}
\toprule
$i$&$v_i=(b,c,d,f,g)$&$\deg v_i$&$i$&$v_i=(b,c,d,f,g)$&$\deg v_i$\\
\midrule
1&$(0,0,0,1,1)$&1&7&$(1,0,1,0,1)$&1\\
2&$(1,1,0,1,0)$&1&8&$(0,1,0,1,1)$&1\\
3&$(1,0,1,1,0)$&2&9&$(0,1,1,0,1)$&1\\
4&$(3,0,1,1,0)$&2&10&$(0,1,1,1,3)$&2\\
5&$(1,0,0,1,0)$&1&11&$(0,1,1,0,0)$&1\\
6&$(1,0,1,0,0)$&1&12&$(0,3,1,1,1)$&2\\
\bottomrule
\end{tabular}
\end{table}

Table~\ref{tab:types} specifies an index set $\sigma_k$ and a subset
$O_k\subseteq\sigma_k$ for each class. Put
\begin{equation}\label{eq:shift}
\beta_k=\sum_{i\in O_k}v_i,\qquad
B_k=\Phi(\beta_k)=\sum_{i\in O_k}V_i.
\end{equation}
The indices in $O_k$ specify the ray coefficients that are positive in
the unshifted cone coordinates.

\begin{table}[H]
\centering
\caption{The twelve normal forms and their boundary shifts.}\label{tab:types}
\begin{tabular}{cclc}
\toprule
$k$&$\sigma_k$&$O_k$&$\deg\beta_k$\\
\midrule
1&$\{1,3,5,6,11\}$&$\{3\}$&2\\
2&$\{1,5,6,7,8\}$&$\{7\}$&1\\
3&$\{1,5,6,8,11\}$&$\varnothing$&0\\
4&$\{1,6,7,8,9\}$&$\{7,9\}$&2\\
5&$\{1,6,8,9,11\}$&$\{9\}$&1\\
6&$\{1,7,8,9,10\}$&$\{10\}$&2\\
7&$\{2,4,5,6,7\}$&$\{4\}$&2\\
8&$\{2,5,6,7,8\}$&$\{2,7\}$&2\\
9&$\{2,5,6,8,11\}$&$\{2\}$&1\\
10&$\{2,6,7,8,9\}$&$\{2,7,9\}$&3\\
11&$\{2,6,8,9,11\}$&$\{2,9\}$&2\\
12&$\{2,8,9,11,12\}$&$\{12\}$&2\\
\bottomrule
\end{tabular}
\end{table}

\begin{theorem}[Normal forms]\label{thm:main}
\begin{enumerate}[label=\textup{(\roman*)},leftmargin=2em]
\item The matrices $V_1,\ldots,V_{12}$ form the unique minimal additive
generating set of $\calS$. Their parameter vectors are exactly the primitive
extreme-ray generators of $C$.
\item There is a disjoint decomposition
\begin{equation}\label{eq:normalforms}
\calS=\bigsqcup_{k=1}^{12}\calF_k,
\qquad
\calF_k=B_k+\sum_{i\in\sigma_k}\N V_i.
\end{equation}
The class and its five coefficients are unique for each $X\in\calS$.
\item The cones $K_k=\cone(v_i:i\in\sigma_k)$ and their faces form a regular
unimodular triangulation of $C$ with respect to $\Z^5$.
\end{enumerate}
\end{theorem}

The unrestricted sum $X=\sum_{i=1}^{12}n_iV_i$ need not be unique. In
\eqref{eq:normalforms}, the allowed generators and the shift are prescribed
by the class.

For $u=(b,c,d,f,g)$ and independent variables $\mathbf y=(y_1,\ldots,y_5)$,
write $\mathbf y^u=y_1^b y_2^c y_3^d y_4^f y_5^g$ and
$w_i=\mathbf y^{v_i}t^{\deg v_i}$. We work with formal series in $t$ whose
coefficients are polynomials in $\mathbf y$.

\begin{corollary}[Generating functions]\label{cor:gf}
The generating function in the five free entries is
\begin{equation}\label{eq:multiGF}
\sum_{u\in C\cap\Z^5}\mathbf y^u t^{\deg u}
=\sum_{k=1}^{12}
\frac{\displaystyle\prod_{i\in O_k}w_i}
{\displaystyle\prod_{i\in\sigma_k}(1-w_i)}.
\end{equation}
In particular, if $a(s)$ counts the matrices of magic sum $2s$, then
\begin{equation}\label{eq:univGF}
F(t):=\sum_{s\ge0}a(s)t^s
=\frac{1+4t+10t^2+4t^3+t^4}{(1-t)^4(1-t^2)}.
\end{equation}
\end{corollary}

To weight the ten upper-triangular entries, replace each $w_i$ on the
right of \eqref{eq:multiGF} by
$t^{\deg v_i}\prod_{1\le r\le q\le4}z_{rq}^{(V_i)_{rq}}$.
The resulting series counts $X\in\calS$ with weight
$t^{\deg X}\prod_{1\le r\le q\le4}z_{rq}^{x_{rq}}$.

\begin{corollary}[Enumeration]\label{cor:count}
For every $s\in\N$,
\begin{equation}\label{eq:count}
a(s)=\frac5{12}s^4+\frac53s^3+\frac{17}{6}s^2+\frac73s
+\frac{7+(-1)^s}{8}.
\end{equation}
For $N\in\N$, let $S_4(N)$ count the squares of magic sum $N$.
Then $S_4(N)=0$ for odd $N$, and for even $N$,
\begin{equation}\label{eq:countN}
S_4(N)=\frac5{192}N^4+\frac5{24}N^3+\frac{17}{24}N^2+\frac76N
+\frac{7+(-1)^{N/2}}8.
\end{equation}
The least quasipolynomial periods are two for $a(s)$ and four for $S_4(N)$.
\end{corollary}

\begin{corollary}[Optimal number of free classes]\label{cor:optimal}
Suppose that $\calS$ is a finite disjoint union
\begin{equation}\label{eq:arbitraryforms}
\calS=\bigsqcup_{j=1}^{r}
\left(A_j+\sum_{\nu=1}^{q_j}\N W_{j\nu}\right),
\end{equation}
where $A_j,W_{j\nu}\in\calS$ and the $q_j$ generators in each class are
linearly independent over $\R$. The ranks $q_j$ may vary between zero and five,
with the empty sum giving a singleton. Then $r\ge12$.
Theorem~\ref{thm:main} attains this bound with $q_j=5$ in every class.
\end{corollary}

\section{The decomposition and its optimality}

\subsection{Coordinates and grading}\label{subsec:coordinates}

Write a general symmetric matrix as
\[
X=\begin{pmatrix}
a&b&c&d\\b&e&f&g\\c&f&h&m\\d&g&m&j
\end{pmatrix}.
\]
The antidiagonal gives $N=2s$, where $s=d+f$. The sum of the row
sums is $8s$, of which the main diagonal contributes $2s$.
Thus the six entries strictly above the diagonal sum to $3s$, giving
$b+c+g+m=2s$. The row equations now give
\begin{equation}\label{eq:entries}
a=s+f-b-c,\quad e=s+d-b-g,\quad
h=b+g-f,\quad m=2s-b-c-g,\quad j=b+c-d.
\end{equation}
For every real $u$, the row, column, and diagonal sums of $\Phi(u)$
in \eqref{eq:phi} are $2s$.
The ten upper-triangular nonnegativity constraints give \eqref{eq:cone}.
The map $\Phi$ and its inverse, which reads the five entries $b,c,d,f,g$,
both have integer coefficients. Hence $\Phi$ identifies $\Z^5$ with the
full lattice of integer solutions to the symmetric magic equations.

The grading $\deg=d+f$ is positive on $C\setminus\{0\}$. If $d+f=0$,
then $d=f=0$, and the upper bounds on $b+c$ and $b+g$ force
$b=c=g=0$. For degree $s$, every matrix entry lies between zero and $2s$,
so there are finitely many lattice points of each degree.
The vector $\eta=(1,1,1,1,1)$ satisfies all defining inequalities strictly.
Thus $C$ is five-dimensional, and it is pointed because it lies in the
nonnegative orthant.

\subsection{Inverse coordinates}\label{subsec:inverse}

We use the linear forms
\begin{equation}\label{eq:comparisons}
Z=b+c+g,\qquad L=b+c-s,\qquad T=b+g-s,\qquad R=Z-s-f.
\end{equation}
Let $M_k$ have columns $v_i$ for $i\in\sigma_k$, in increasing order.
Table~\ref{tab:inverse} defines the linear map $\lambda^{(k)}$.
Its coordinate $\lambda_i^{(k)}(u)$ is indexed by the corresponding
column $v_i$.

\begin{table}[htbp]
\centering
\caption{Integer inverse coordinates. The letters $a,e,h,m,j$ denote
the linear forms in \eqref{eq:entries}.}\label{tab:inverse}
\small
\begin{tabular}{clc}
\toprule
$k$&$\lambda^{(k)}(u)$&$\det M_k$\\
\midrule
1&$(g,\ s-Z,\ j,\ h,\ c)$&$1$\\
2&$(-L,\ b-d,\ e,\ T,\ c)$&$-1$\\
3&$(-L,\ f-g,\ h,\ Z-s,\ -T)$&$-1$\\
4&$(-L,\ s-g,\ T,\ j,\ d-b)$&$-1$\\
5&$(-L,\ b,\ j,\ g-f,\ -T)$&$1$\\
6&$(m,\ b,\ j,\ e,\ g-s)$&$-1$\\
7&$(c,\ b-s,\ a,\ e,\ g)$&$1$\\
8&$(L,\ f-c,\ e,\ T,\ s-b)$&$1$\\
9&$(L,\ -R,\ h,\ g,\ -T)$&$1$\\
10&$(L,\ m,\ T,\ a,\ c-f)$&$1$\\
11&$(L,\ s-c,\ a,\ R,\ -T)$&$-1$\\
12&$(b,\ a,\ h,\ m,\ c-s)$&$-1$\\
\bottomrule
\end{tabular}
\end{table}

The table satisfies $M_k\lambda^{(k)}(u)=u$ for all $u\in\R^5$.
For example,
row~1 gives
\begin{equation}\label{eq:outerexample}
u=g v_1+(s-Z)v_3+jv_5+hv_6+cv_{11}.
\end{equation}
The first, second, and fifth entries on the right are $j+h+(s-Z)$,
$c$, and $g$. By \eqref{eq:entries}, the first of these is $b$.
The remaining entries are $(s-Z)+h+c=d$ and $g+(s-Z)+j=f$.
Thus $M_k$ is invertible with $\lambda^{(k)}(u)=M_k^{-1}u$.
Both maps have integer coefficients, so $|\det M_k|=1$.

Substitution in \eqref{eq:cone} gives $v_i\in C$ for every $i$, hence $K_k\subseteq C$.
Define a subset of $K_k$ by
\begin{equation}\label{eq:halfopen}
K_k^\circledast=
\left\{u\in K_k:\lambda_i^{(k)}(u)>0\text{ for all }i\in O_k\right\}.
\end{equation}
The set $K_3^\circledast$ equals the closed cone $K_3$ because $O_3$ is
empty. At a lattice point, the strict inequalities require
$\lambda_i^{(k)}(u)\ge1$ for $i\in O_k$. Therefore
\begin{equation}\label{eq:latticeclass}
K_k^\circledast\cap\Z^5
=\beta_k+\sum_{i\in\sigma_k}\N v_i.
\end{equation}
We next show that the twelve sets $K_k^\circledast$ partition $C$.

\subsection{Four outer regions and eight central regions}\label{subsec:partition}

For a point of $C$, at most one of the following inequalities holds.
\begin{equation}\label{eq:outerconditions}
Z<s,\qquad g>s,\qquad b>s,\qquad c>s.
\end{equation}
The first excludes the other three since $b,c,g\ge0$. Any two of the
remaining inequalities would force one of $b+c$, $b+g$, or $Z$ above
$2s$, contrary to \eqref{eq:cone}. Rows $1,6,7,12$ of
Table~\ref{tab:inverse} identify these regions with
$K_1^\circledast,K_6^\circledast,K_7^\circledast,K_{12}^\circledast$,
respectively. In each row the coefficient indexed by $O_k$ must
be positive, and the other four are entries of $X$.

All remaining points lie in the central region
\begin{equation}\label{eq:central}
C_0=\{u\in C:b\le s,\ c\le s,\ g\le s,\ Z\ge s\}.
\end{equation}
The remaining cones all lie in $C_0$ because their generators satisfy
\eqref{eq:central}. For a real number $x$, write $x_+=\max(x,0)$.
An integral equivalence is a linear isomorphism whose matrix and inverse
have integer entries.

\begin{lemma}[A common cone for the central regions]\label{lem:central}
For each $\varepsilon,\tau\in\{-1,1\}$, the cone
$C_0\cap\{\varepsilon L\ge0,\ \tau T\ge0\}$ is integrally equivalent to
$\R_{\ge0}^{2}\times W$, where
\[
W=\{(p,q,r)\in\R^3:p,q\ge0,\ 0\le r\le p+q\}.
\]
The integral linear coordinates are
\begin{equation}\label{eq:centralcoords}
\alpha=|L|,\quad \beta=|T|,\quad
p=f-|L|,\quad q=d-|T|,\quad r=b-L_+-T_+.
\end{equation}
The sets $W\cap\{r\ge q\}$ and $W\cap\{r<q\}$ have respective free
coordinates
\begin{equation}\label{eq:centralfree}
(p+q-r,\ r-q,\ q),\qquad (p,\ q-r,\ r).
\end{equation}
The first tuple ranges over $\R_{\ge0}^3$. The second ranges over the
triples in $\R_{\ge0}^3$ with strictly positive second coordinate.
\end{lemma}
\begin{proof}
The four diagonal-entry inequalities $f\pm L\ge0$ and $d\pm T\ge0$
are equivalent to $p,q\ge0$ and imply $d,f\ge0$.
Substitute $c=s+L-b$, $g=s+T-b$, and $Z=2s+L+T-b$ into
the remaining inequalities. The conditions $b,c,g\ge0$ and $Z\le2s$,
together with $b,c,g\le s$ and $Z\ge s$, give exactly
\[
\max(0,L,T,L+T)\le b\le
\min(s,s+L,s+T,s+L+T).
\]
The lower bound is $L_++T_+$. The interval has length
$s-|L|-|T|=p+q$, so these inequalities are equivalent to $0\le r\le p+q$.

After fixing the signs of $L$ and $T$, the inverse map is integral:
\[
\begin{gathered}
d=q+\beta,\quad f=p+\alpha,\quad s=d+f,\quad
b=r+L_++T_+,\\
c=s+L-b,\qquad g=s+T-b,
\end{gathered}
\]
where $L$ and $T$ are the chosen signed values of $\alpha$ and $\beta$.
For each fixed sign choice, the forward and inverse formulas extend
to integer linear maps on the ambient spaces. They are mutual inverses
and identify the two ambient lattices.
For the first tuple $(x,y,z)$ in \eqref{eq:centralfree}, its inverse is
$(p,q,r)=(x+y,z,y+z)$. For the second tuple it is
$(p,q,r)=(x,y+z,z)$. In the second case $r<q$ is equivalent to $y>0$.
\end{proof}

Assign the sign boundaries to $L\le0$ and $T\le0$. In each of the four
sign regions, \eqref{eq:centralcoords} gives
\begin{equation}\label{eq:centralbridge}
r-q=\begin{cases}
b-d,&L\le0,\ T>0,\\
f-g,&L\le0,\ T\le0,\\
f-c,&L>0,\ T>0,\\
-R,&L>0,\ T\le0.
\end{cases}
\end{equation}
These identities give the last comparison in Table~\ref{tab:decision}.

\begin{table}[H]
\centering
\caption{A disjoint decision table for the central region $C_0$.}\label{tab:decision}
\begin{tabular}{cccc}
\toprule
First comparison&Second comparison&Final comparison&Class $k$\\
\midrule
$L\le0$&$T>0$&$b\ge d$&2\\
$L\le0$&$T>0$&$b<d$&4\\
$L\le0$&$T\le0$&$g\le f$&3\\
$L\le0$&$T\le0$&$g>f$&5\\
$L>0$&$T>0$&$c\le f$&8\\
$L>0$&$T>0$&$c>f$&10\\
$L>0$&$T\le0$&$R\le0$&9\\
$L>0$&$T\le0$&$R>0$&11\\
\bottomrule
\end{tabular}
\end{table}

Let $(x,y,z)$ denote the corresponding tuple from \eqref{eq:centralfree}.
Together with $\alpha,\beta$, it gives the inverse rows in the following
column orders.
\[
\begin{array}{c|c@{\qquad}c|c}
k& M_k^{-1}u&k&M_k^{-1}u\\ \hline
2,8&(\alpha,y,z,\beta,x)&4,10&(\alpha,z,\beta,x,y)\\
3,9&(\alpha,y,z,x,\beta)&5,11&(\alpha,z,x,y,\beta).
\end{array}
\]
For example, if $L\le0$ and $T>0$, then $\alpha=-L$, $\beta=T$,
and $r-q=b-d$. The first tuple is $(c,b-d,e)$, which gives row~2.
The second tuple is $(j,d-b,s-g)$, which gives row~4. In each sign
region, $\alpha>0$ is required when $L>0$ and $\beta>0$ when $T>0$.
The second tuple also requires $y>0$. Their positions in the displayed
column orders give exactly the sets $O_k$ in Table~\ref{tab:types}.
Conversely, the signed $L,T$ and $r-q$ or $q-r$ entries in each row
of Table~\ref{tab:inverse}, with the strict positions $O_k$, impose
the corresponding three inequalities in Table~\ref{tab:decision}.
Thus the table
assigns each point of $C_0$ to exactly one $K_k^\circledast$.

This proves the exact real partition
\begin{equation}\label{eq:realpartition}
C=\bigsqcup_{k=1}^{12}K_k^\circledast.
\end{equation}
Intersecting with $\Z^5$ and using \eqref{eq:latticeclass} proves
Theorem~\ref{thm:main}(ii). For a given $X\in\calS$, read
\[
u=(x_{12},x_{13},x_{14},x_{23},x_{24}).
\]
Compute $s,Z,L,T,R$,
select an outer class by \eqref{eq:outerconditions} or a central class by
Table~\ref{tab:decision}, and evaluate the corresponding row of
Table~\ref{tab:inverse}. Subtract one from $\lambda_i^{(k)}(u)$ for each
$i\in O_k$ to obtain the five parameters in \eqref{eq:normalforms}.

\begin{example}\label{ex:normalform}
Take $u=(2,3,2,2,1)$. Then $s=4$, $Z=6$, $L=1$, $T=-1$, and $R=0$.
The point belongs to $K_9^\circledast$. Its inverse coordinates in
the order $(2,5,6,8,11)$ are $(1,0,1,1,1)$. Thus
\[
\Phi(u)=V_2+V_6+V_8+V_{11}
=\begin{pmatrix}
1&2&3&2\\2&3&2&1\\3&2&1&2\\2&1&2&3
\end{pmatrix}.
\]
Since $B_9=V_2$, the five parameters of the matrix in class~9 are
$(0,0,1,1,1)$. Every line sum is eight.
\end{example}

\subsection{Extreme rays and the Hilbert basis}\label{subsec:hilbert}

The cone $C$ is defined by the nonnegativity of the following ten forms.
\begin{equation}\label{eq:entryforms}
\begin{aligned}
\ell_1&=b,&\ell_2&=c,&\ell_3&=d,&\ell_4&=f,&\ell_5&=g,\\
\ell_6&=b+c-d,&\ell_7&=b+g-f,&
\ell_8&=d+2f-b-c,\\
\ell_9&=2d+f-b-g,&
\ell_{10}&=2d+2f-b-c-g.
\end{aligned}
\end{equation}
Table~\ref{tab:active} lists, for each $v_i$, four forms whose coefficient
rows are linearly independent and vanish on $v_i$. Their common kernel is the line spanned by $v_i$ and
intersects $C$ in the ray $\R_{\ge0}v_i$ because $C$ is pointed and
$v_i\in C$. The intersection is a face: it is the zero set in $C$ of the
sum of the four nonnegative forms. Each $v_i$ has a coordinate equal to
one, so it is primitive.

\begin{table}[H]
\centering
\caption{Independent active constraints for the extreme rays.}\label{tab:active}
\begin{tabular}{cc|cc|cc}
\toprule
$i$&Active indices&$i$&Active indices&$i$&Active indices\\
\midrule
1&$1,2,3,7$&5&$2,3,5,7$&9&$1,4,6,10$\\
2&$3,5,7,8$&6&$2,4,5,6$&10&$1,6,9,10$\\
3&$2,5,6,7$&7&$2,4,6,9$&11&$1,4,5,6$\\
4&$2,5,8,9$&8&$1,3,7,10$&12&$1,7,8,10$\\
\bottomrule
\end{tabular}
\end{table}

The partition gives $C=\cone(v_1,\ldots,v_{12})$ and
$C\cap\Z^5=\sum_i\N v_i$. Thus the list of extreme rays is complete.
If $v_i=u+w$ for nonzero $u,w\in C\cap\Z^5$, extremality puts both
summands on the same ray. Primitivity makes them positive integer
multiples of $v_i$, whose sum cannot equal $v_i$. Every other nonzero
lattice point is a sum of at least two generators. The indecomposable
nonzero elements must occur in every additive generating set, so
\[
\Hilb(C)=\{v_1,\ldots,v_{12}\}.
\]
This proves Theorem~\ref{thm:main}(i).

\subsection{Generating functions and coefficient extraction}\label{subsec:countproof}

Summing over the five parameters in each normal form gives
\eqref{eq:multiGF}. In the central region $s=\alpha+\beta+p+q$.
For each of $L,T$, a positive absolute value has two sign choices and
zero has one, contributing $(1+t)/(1-t)$. The two tuples in
\eqref{eq:centralfree} give $p+q=x+y+z$, with $y>0$ in the second piece,
so $W$ contributes $(1+t)/(1-t)^3$. Thus the central contribution is
$(1+t)^3/(1-t)^5$. Each of the four outer classes has four degree-one
generators, one degree-two generator, and a degree-two shift. Therefore
\begin{align}
F(t)
&=\frac{1+3t+3t^2+t^3}{(1-t)^5}
 +\frac{4t^2}{(1-t)^4(1-t^2)}\notag\\
&=\frac{(1+t)^4+4t^2}{(1-t)^4(1-t^2)},\label{eq:groupedGF}
\end{align}
which proves Corollary~\ref{cor:gf}. Expanding into partial fractions gives
\begin{equation}\label{eq:partial}
\begin{split}
F(t)={}&\frac{10}{(1-t)^5}-\frac{15}{(1-t)^4}
+\frac{13}{2(1-t)^3}-\frac{3}{4(1-t)^2}\\
&+\frac1{8(1-t)}+\frac1{8(1+t)}.
\end{split}
\end{equation}
Using $[t^s](1-t)^{-r}=\binom{s+r-1}{r-1}$ gives
\eqref{eq:count}. Taking $s=N/2$ gives \eqref{eq:countN} for even $N$.
Odd magic sums are excluded by the antidiagonal condition.
For a quasipolynomial of degree at most four and period $p$, summing
over residue classes gives a denominator dividing $(1-z^p)^5$.
The numerator in \eqref{eq:univGF} is $4$ at $t=-1$, so $F(t)$ has
a pole there. The magic-sum series is $F(z^2)$, which has a pole at
$z=\mathrm i$. Thus the periods of $a(s)$ and $S_4(N)$ are multiples
of two and four, respectively. The formulas attain these periods.

The first values of $a(s)$, starting at $s=0$, are
$1,8,37,112,269,552,1017,1728,2761$.
Summing the number of admissible values of $g$ in \eqref{eq:cone} gives
\begin{equation}\label{eq:directcount}
\begin{split}
a(s)=\sum_{d=0}^{s}\sum_{b=0}^{2s}
\sum_{c=\max(0,d-b)}^{s+f-b}
\Bigl[\min(s+d-b,2s-b-c)-\max(0,f-b)+1\Bigr]_+,\\[-2pt]
\hfill f=s-d,
\end{split}
\end{equation}
where $[r]_+=\max(r,0)$ and an empty range contributes zero.

\paragraph{Relation to the semimagic series.}
Consider a symmetric semimagic square with nonnegative entries. We use
the two disjoint perfect matchings $12,34$ and $13,24$
(cf.~\cite[Corollary~4.10]{Stanley76}). In each pair, subtract the smaller
entry from both entries and from their transposed positions. Each
subtraction removes a multiple of a perfect-matching matrix of row
sum one. The two supports are disjoint, so the original square is
recovered uniquely from the two removed amounts and the residual
square. Conversely, the two amounts can be chosen independently
for any residual square with at most one nonzero entry in each pair.
They contribute $(1-t)^{-2}$ to the row-sum series.

Choose a surviving edge from each pair, using a fixed edge when both
entries in its pair are zero. Write $b',c'$ for their weights. The
chosen edges meet at a vertex. Let $a'$ be the diagonal entry there
and let $d',f'$ be the weights of the third opposite pair, with $d'$
incident to that vertex. If $n$ is the row sum of the residual square,
the diagonal entries at the other ends of the edges with weights
$b',c',d'$ are $n-b'-f'$, $n-c'-f'$, and $n-d'$, respectively.
The row at the common vertex determines $n$. Nonnegativity gives
\[
n=a'+b'+c'+d',
\qquad 0\le f'\le a'+d'+\min(b',c'),
\]
where $a',b',c',d'$ are nonnegative. Each positive $b'$ or $c'$ has two
possible positions, and a zero weight has one. When $f'\le a'+d'$,
summing $(a'+d'+1)t^{a'+d'}$ over $a',d'$ gives $(1+t)/(1-t)^3$.
The factors for $b',c'$ are both $(1+t)/(1-t)$.
When $f'>a'+d'$, put $k=f'-a'-d'\ge1$. Then $b',c'\ge k$, so there
are four choices of positions. The remaining parameters
$a',d',b'-k,c'-k$ range freely over $\N$, and their sum plus $2k$ is
$n$. This contributes $4t^2/((1-t)^4(1-t^2))$.
The two contributions are precisely those in \eqref{eq:groupedGF}.
Consequently $G(t)=(1-t)^{-2}F(t)$. In particular, this explains the
common numerator in our series and the semimagic series of
Carlitz \cite{Carlitz66,Stanley76}.

\subsection{Why twelve free classes are necessary}\label{subsec:optimal}

\begin{proof}[Proof of Corollary~\ref{cor:optimal}]
Put $e_{j\nu}=\deg W_{j\nu}>0$. Disjointness and linear independence
give
\[
F(t)=\sum_{j=1}^{r}
       \frac{t^{\deg A_j}}{\prod_{\nu=1}^{q_j}(1-t^{e_{j\nu}})}.
\]
Let $J_5=\{j:q_j=5\}$, $r_5=|J_5|$, and
$D_j=\prod_{\nu=1}^5 e_{j\nu}$ for $j\in J_5$. Multiplying by
$(1-t)^5$ and letting $t$ tend to one from below removes the
lower-rank terms. By \eqref{eq:univGF},
\begin{equation}\label{eq:degreebudget}
\sum_{j\in J_5}\frac1{D_j}=10,\qquad
r_5-10=\sum_{j\in J_5}\left(1-\frac1{D_j}\right).
\end{equation}
In particular, $r_5\ge10$.

Each extreme ray of $C$ must contain a generator $\Phi^{-1}(W_{j\nu})$.
To see this, take a positive integer $n$ with
$n\deg v_i>\max_j\deg A_j$. The representation of $nV_i$ in
\eqref{eq:arbitraryforms} uses at least one generator with positive
coefficient. Applying $\Phi^{-1}$, every nonzero summand lies on the
ray of $v_i$, since it is a face of $C$. On each of the four rays
of $v_3,v_4,v_{10},v_{12}$ a nonzero lattice vector is a positive
integer multiple of the primitive generator. Each of those rays
therefore requires a generator of even degree in the decomposition.

Suppose $r\le11$. If every class has rank five, the second sum in
\eqref{eq:degreebudget} is at most one. There are at most two indices
with $D_j>1$, since each contributes at least $1/2$. There must be
at least one to account for the four rays. A single contribution
lies strictly between zero and one, contrary to the integrality of
$r_5-10$. If there are two, each equals $1/2$. Both products then
equal two, and the two classes have only one even-degree generator
each. This again falls short of four.

If some class has rank below five, then $r_5=10$ and there is just
one lower-rank class. Every full-rank class has $D_j=1$ by
\eqref{eq:degreebudget}. The four even-degree generators required
on the distinct rays must all belong to the lower-rank class.
That class has rank four, and its generating-function summand has
a pole of order four at $t=-1$. Its monomial numerator is nonzero
there, and the full-rank summands are analytic there because all
their generator degrees are one. But \eqref{eq:univGF} has a
simple pole at $-1$. This contradiction finishes the proof.
\end{proof}

\section{Geometry and consequences}

\subsection{Regularity of the triangulation}\label{subsec:regular}

We realize the triangulation by projecting the lower facets of a lifted
polytope, as in \cite{DRS10}. Set $\delta_i=\deg v_i$ and $p_i=v_i/\delta_i$,
so $P=\conv(p_1,\ldots,p_{12})$. Choose
\begin{equation}\label{eq:heights}
(H_1,\ldots,H_{12})=(0,-1,6,0,0,-2,-3,-2,-1,4,2,6)
\end{equation}
and lift $p_i$ to height $H_i/\delta_i$. Write
\[
\widehat P=\conv\{(p_i,H_i/\delta_i):1\le i\le12\}.
\]
For each $k$, let $h_k$ be the linear form whose coefficients in the order
$(b,c,d,f,g)$ are given in Table~\ref{tab:lifting-perturbation}. Substitution gives
\begin{equation}\label{eq:slacks}
H_i-h_k(v_i)=0\quad(i\in\sigma_k),\qquad
H_i-h_k(v_i)>0\quad(i\notin\sigma_k).
\end{equation}
For example, the slacks for $k=1$, in ray order, are
$(0,3,0,10,0,0,7,2,5,18,0,12)$.

\begin{table}[htbp]
\centering
\caption{Lifting forms $h_k$ and perturbation coordinates $M_k^{-1}\rho$.}
\label{tab:lifting-perturbation}
\small
\begin{tabular}{ccc}
\toprule
$k$&Coefficients of $h_k$&$M_k^{-1}\rho$\\
\midrule
1&$(-8,-4,6,8,-8)$&$(2,-1,3,1,2)$\\
2&$(-1,-2,-1,1,-1)$&$(1,1,3,-1,2)$\\
3&$(-6,-2,4,6,-6)$&$(1,2,1,1,1)$\\
4&$(-4,-2,2,1,-1)$&$(1,4,-1,3,-1)$\\
5&$(-6,-2,4,3,-3)$&$(1,3,3,-2,1)$\\
6&$(-4,-2,-4,-5,5)$&$(5,3,3,3,-4)$\\
7&$(2,-1,-4,-2,-1)$&$(2,-3,5,3,2)$\\
8&$(0,-1,-2,0,-1)$&$(-1,2,3,-1,3)$\\
9&$(-5,-1,3,5,-6)$&$(-1,3,1,2,1)$\\
10&$(0,2,-2,-3,-1)$&$(-1,5,-1,5,-2)$\\
11&$(-2,2,0,-1,-3)$&$(-1,4,5,-3,1)$\\
12&$(-2,6,-4,-5,-3)$&$(3,5,1,5,-4)$\\
\bottomrule
\end{tabular}
\end{table}

Dividing \eqref{eq:slacks} by $\delta_i>0$ gives the supporting hyperplane
$y=h_k(x)$ of $\widehat P$ and the face
\begin{equation}\label{eq:lowerface}
\widehat P\cap\{(x,y):y=h_k(x)\}
=\conv\{(p_i,H_i/\delta_i):i\in\sigma_k\}.
\end{equation}
The nonsingularity of $M_k$ makes its five normalized columns affinely
independent. The face in \eqref{eq:lowerface} therefore has dimension four.
All the other lifted points lie strictly above its hyperplane, so
$\widehat P$ has dimension five in $\{\deg x=1\}\times\R$, and this face is
a lower facet.

Its projection is $K_k\cap P$. These simplices cover $P$ by
\eqref{eq:realpartition}. The lower envelope of $\widehat P$ is the graph of the minimum height
above each $x\in P$. Projection is injective on this graph, so the
intersection of two projected lower facets is the projection of their
common face. Their relative interiors are therefore disjoint. The twelve facets in \eqref{eq:lowerface}
thus account for all the lower facets. Coning over their projections gives
the desired regular triangulation, and the determinants in
Table~\ref{tab:inverse} establish unimodularity in $\Z^5$. This proves
Theorem~\ref{thm:main}(iii).

\subsection{Boundary assignment}\label{subsec:boundary}

We now recover the sets $O_k$ from a single interior perturbation.
We use the half-open subdivision rule of \cite[Section~2.2]{KV08}
(see also \cite{BS07}). Take
\begin{equation}\label{eq:rho}
\rho=(3,2,2,4,2),\qquad
(\ell_1(\rho),\ldots,\ell_{10}(\rho))=(3,2,2,4,2,3,1,5,3,5).
\end{equation}
All coordinates in the second tuple are positive. Thus $u+\varepsilon\rho$
lies in the interior of $C$ for every $u\in C$ and $\varepsilon>0$.
As Table~\ref{tab:lifting-perturbation} shows, $\lambda_i^{(k)}(\rho)\ne0$ for every
$k$ and $i\in\sigma_k$. Fix $u\in C$. For all sufficiently small
$\varepsilon>0$, the sign of
\[
\lambda_i^{(k)}(u+\varepsilon\rho)
=\lambda_i^{(k)}(u)+\varepsilon\lambda_i^{(k)}(\rho)
\]
is the sign of $\lambda_i^{(k)}(u)$ when this is nonzero, and the sign of
$\lambda_i^{(k)}(\rho)$ when $\lambda_i^{(k)}(u)=0$. There are finitely many
coordinates, so the perturbed point avoids all subdivision walls and lies
in the interior of exactly one $K_k$. It lies in that interior precisely when
\[
\lambda_i^{(k)}(u)\ge0\quad(i\in\sigma_k),\qquad
\lambda_i^{(k)}(u)>0\quad
\bigl(i\in\sigma_k,\ \lambda_i^{(k)}(\rho)<0\bigr).
\]
The negative coordinates in Table~\ref{tab:lifting-perturbation} are indexed by $O_k$.
Thus the perturbation gives \eqref{eq:halfopen} for all points of $C$,
including the extreme rays and the origin. For $u\in C\cap\Z^5$, one may take
$\varepsilon=1/6$ uniformly: the unperturbed coordinates are integers,
and Table~\ref{tab:lifting-perturbation} bounds the absolute perturbation
coordinates by five.

Different parts of a shared face can belong to different half-open cones.
For example, $K_2\cap K_8=\cone(v_5,v_6,v_7,v_8)$. The point
$v_5+v_6+v_7+v_8=(3,1,2,2,2)$ in its relative interior belongs to
$K_2^\circledast$. The ray generator $v_5=(1,0,0,1,0)$ on its boundary
belongs to $K_3^\circledast$.

\subsection{Fundamental parallelepipeds and the shifts}\label{subsec:remainders}

For a simplicial cone $K=\cone(r_1,\ldots,r_5)$ with integral generators,
write
\[
\Pi_K=\left\{\sum_{j=1}^5\theta_jr_j:0\le\theta_j<1\right\}.
\]
For each $K_k$, unimodularity gives $\Pi_{K_k}\cap\Z^5=\{0\}$. To match
\eqref{eq:halfopen}, use $(0,1]$ in place of $[0,1)$ for coordinates
indexed by $O_k$. The resulting parallelepiped contains the single lattice
point $\beta_k=\sum_{i\in O_k}v_i$, which accounts for the numerator of
the $k$th summand in \eqref{eq:multiGF}.

A lattice point in a fundamental parallelepiped may decompose in the
full monoid. For instance, take
\[
K'=\cone(v_1,v_7,v_8,v_{10},v_{11}),\qquad
|\det(v_1,v_7,v_8,v_{10},v_{11})|=2.
\]
The point
\begin{equation}\label{eq:localexample}
p=(0,1,1,1,2)=\tfrac12(v_1+v_{10}+v_{11})=v_1+v_9
\end{equation}
belongs to $\Pi_{K'}$ and decomposes as a sum of two nonzero points of
$C\cap\Z^5$.

\subsection{Comparison with earlier counts}\label{subsec:comparison}

Both preprint versions \cite{BCCGv2,BCCGv3} impose the same symmetry
and diagonal conditions as above. Their Section~4 formulas for $S_4(N)$
are on page~8 of v2 and page~9 of v3. The values also appear in
OEIS A093198, indexed by $s=N/2$
\cite{OEISA093198}. Table~\ref{tab:comparison} compares them with the
values from Corollary~\ref{cor:count}.

\begin{table}[htbp]
\centering
\caption{Counts at four small magic sums.}\label{tab:comparison}
\begin{tabular}{@{}lrrrr@{}}
\toprule
Magic sum $N$&2&4&6&8\\
\midrule
Formula~\eqref{eq:countN}&8&37&112&269\\
arXiv versions v2 and v3, and OEIS A093198&11&53&164&397\\
\bottomrule
\end{tabular}
\end{table}

At $N=2$, the matrix equations give the count directly.
Use the entries $b,c,d,f,g,m$ in the positions $(12,13,14,23,24,34)$,
respectively. The antidiagonal condition gives $d+f=1$. The total of the
four row sums is eight, and the trace is two, so the six entries above
the main diagonal sum to three. Hence $b+c+g+m=2$.

If $(d,f)=(0,1)$, nonnegativity of the second and third diagonal entries
requires $b+g\le1$ and $c+m\le1$. Both are therefore equalities.
There are four choices, and the first and fourth diagonal entries are
nonnegative for each choice. If $(d,f)=(1,0)$, the first and fourth row
equations give $b+c\le1$ and $g+m\le1$. Both are equalities, giving four
more choices with nonnegative diagonal entries in the other two rows.
The row sums determine the diagonal entries uniquely. For each choice
the trace is $8-2\cdot3=2$, and the antidiagonal sum is $2(d+f)=2$.
This gives eight symmetric magic squares of magic sum two.

\subsection{Positive entries and a reflexive polytope}\label{subsec:reflexive}

The all-ones matrix $J=\Phi(\eta)$ has magic sum four, where
$\eta=(1,1,1,1,1)$. For $N\ge4$, subtracting $J$ gives a bijection from
symmetric magic squares with positive entries and magic sum $N$ to those
with nonnegative entries and magic sum $N-4$. For $N<4$ there are
no squares with positive entries. In parameter coordinates,
the bijection is expressed by
\begin{equation}\label{eq:interiorshift}
\intr(C)\cap\Z^5=\eta+(C\cap\Z^5).
\end{equation}
The interior of $C$ consists of the points where all the forms in
\eqref{eq:entryforms} are positive. In fact, the forms are nonzero,
so equality in any defining inequality puts a point on the boundary.
Strict inequalities persist in a neighborhood. Each form is integer-valued
on $\Z^5$ and satisfies $\ell_i(\eta)=1$.
For $u\in\Z^5$, subtracting $\eta$ turns the conditions $\ell_i(u)>0$
into $\ell_i(u-\eta)\ge0$.

Let $\Lambda=\ker(\deg)\cap\Z^5$, $\Lambda_\R=\Lambda\otimes\R$, and
$\Lambda^*=\operatorname{Hom}(\Lambda,\Z)$. A full-dimensional lattice
polytope $Q\subset\Lambda_\R$ containing the origin in its interior is
reflexive if its polar
\[
Q^\vee=\{y\in\Lambda^*\otimes\R:\langle y,z\rangle\le1
\text{ for every }z\in Q\}
\]
is a lattice polytope in the dual lattice. In dimension four the Ehrhart
numerator $h_Q^*(t)$ is defined by
\begin{equation}\label{eq:hstardef}
\sum_{n\ge0}\#(nQ\cap\Lambda)t^n=\frac{h_Q^*(t)}{(1-t)^5}.
\end{equation}

\begin{proposition}\label{prop:reflexive}
The polytope $Q=2P-\eta$ is reflexive in $\Lambda$ and has
\begin{equation}\label{eq:hstar}
h_Q^*(t)=1+32t+94t^2+32t^3+t^4.
\end{equation}
\end{proposition}
\begin{proof}
Every vertex of $2P$ is $2v_i/\delta_i$ for some $\delta_i\in\{1,2\}$.
Thus $Q$ is a four-dimensional lattice polytope in $\Lambda$. Since
$\deg\eta=2$ and $\ell_i(\eta)=1$ for every $i$,
\[
Q=\{z\in\Lambda_\R:\ell_i(z)\ge-1\text{ for }1\le i\le10\}.
\]
All inequalities are strict at the origin. To take the polar, write
the displayed system as $-\ell_i(z)\le1$. The polar of the convex
hull of the covectors $-\ell_i|_\Lambda$ is $Q$. This hull contains
the origin, since $\ell_3+\ell_4=d+f$ vanishes on $\Lambda$.
Bipolarity therefore gives
\begin{equation}\label{eq:polar}
Q^\vee=\conv\{-\ell_i|_\Lambda:1\le i\le10\}.
\end{equation}
The listed covectors belong to $\Lambda^*$, so $Q$ is reflexive.

Translation by $n\eta$ identifies $nQ\cap\Lambda$ with the degree-$2n$
lattice points of $C$. Formula~\eqref{eq:count} yields
\[
\#(nQ\cap\Lambda)=a(2n)
=\frac{20}{3}n^4+\frac{40}{3}n^3+\frac{34}{3}n^2+\frac{14}{3}n+1.
\]
The coefficient of $t^j$ in $(1-t)^5\sum_{n\ge0}a(2n)t^n$ is
$\sum_{k=0}^{\min(5,j)}(-1)^k\binom{5}{k}a(2j-2k)$.
For $j\ge5$ this is a fifth difference of a polynomial of degree
four and vanishes. Substitution for $0\le j\le4$ gives
\eqref{eq:hstar}.
\end{proof}

The Ehrhart numerator of $Q$ in \eqref{eq:hstar} is palindromic and
unimodal. The series in \eqref{eq:univGF} is the Ehrhart series of the rational
section $P$.
By \eqref{eq:interiorshift}, the series for positive entries is $t^2F(t)$.
The identity $F(t^{-1})=-t^2F(t)$ is the corresponding five-dimensional
cone reciprocity relation. Davis \cite{Davis15} studied Ehrhart series
of polytopes of symmetric doubly stochastic matrices with only row and
column sums prescribed.

\appendix
\section{The twelve basic matrices}\label{app:matrices}

These are the matrices $V_i=\Phi(v_i)$ for the vectors in Table~\ref{tab:rays}.
\begin{align*}
V_1&=\begin{pmatrix}2&0&0&0\\0&0&1&1\\0&1&0&1\\0&1&1&0\end{pmatrix},&
V_2&=\begin{pmatrix}0&1&1&0\\1&0&1&0\\1&1&0&0\\0&0&0&2\end{pmatrix},&
V_3&=\begin{pmatrix}2&1&0&1\\1&2&1&0\\0&1&0&3\\1&0&3&0\end{pmatrix},\\[2pt]
V_4&=\begin{pmatrix}0&3&0&1\\3&0&1&0\\0&1&2&1\\1&0&1&2\end{pmatrix},&
V_5&=\begin{pmatrix}1&1&0&0\\1&0&1&0\\0&1&0&1\\0&0&1&1\end{pmatrix},&
V_6&=\begin{pmatrix}0&1&0&1\\1&1&0&0\\0&0&1&1\\1&0&1&0\end{pmatrix},\\[2pt]
V_7&=\begin{pmatrix}0&1&0&1\\1&0&0&1\\0&0&2&0\\1&1&0&0\end{pmatrix},&
V_8&=\begin{pmatrix}1&0&1&0\\0&0&1&1\\1&1&0&0\\0&1&0&1\end{pmatrix},&
V_9&=\begin{pmatrix}0&0&1&1\\0&1&0&1\\1&0&1&0\\1&1&0&0\end{pmatrix},\\[2pt]
V_{10}&=\begin{pmatrix}2&0&1&1\\0&0&1&3\\1&1&2&0\\1&3&0&0\end{pmatrix},&
V_{11}&=\begin{pmatrix}0&0&1&1\\0&2&0&0\\1&0&0&1\\1&0&1&0\end{pmatrix},&
V_{12}&=\begin{pmatrix}0&0&3&1\\0&2&1&1\\3&1&0&0\\1&1&0&2\end{pmatrix}.
\end{align*}

\end{document}